\documentclass[11pt]{article}%
\usepackage{amsmath}
\usepackage{graphicx}
\usepackage{amsfonts}
\usepackage{amssymb}
\usepackage{a4wide}%
\newtheorem{theorem}{Theorem}

\newtheorem{corollary}{Corollary}

\newtheorem{definition}{Definition}
\newtheorem{example}{Example}

\newtheorem{lemma}{Lemma}

\newtheorem{proposition}{Proposition}
\newtheorem{remark}{Remark}

\begin{document}

\title{Strongly adequate functions on Banach spaces}
\author{Michel Volle\thanks{Department of Mathematics, University of Avignon, France,
e-mail: \texttt{michel.volle@univ-avignon.fr}.}~ and Constantin
Z\u{a}linescu\thanks{Faculty of Mathematics, University Al. I. Cuza, Ia\c{s}i,
Romania, e-mail: \texttt{zalinesc@uaic.ro}.}}
\maketitle

\textbf{Abstract} The notion of adequate function has been recently introduced
in order to characterize the essentially strictly functions on a reflexive
Banach space among the weakly lower semicontinuous ones. In this paper we
reinforce this concept and show that a lower semicontinuous function is
essentially firmly subdifferentiable if and only if it is strongly adequate.

\bigskip


\bigskip


\section{Introduction}

The notion of adequate function on a reflexive Banach space has been recently
introduced to obtain a characterization for the class of essentially strictly
convex functions (in the sense of \cite{R3}) among the weakly lower
semicontinuous ones (\cite[Th.\ 1]{R16}). In the present paper we reinforce
this concept (Definition \ref{def1}) in order to treat lower semicontinuous
(lsc) functions on general Banach spaces instead of weakly lower
semicontinuous ones. The corresponding concept of convexity no longer produces
the class of essentially strictly convex functions but the class of
essentially firmly subdifferentiable (convex) functions that we introduce for
this purpose (Definition \ref{def3}); this notion benefits from nice
properties in terms of optimization problems (Proposition \ref{prop5b}). We
prove in Theorem \ref{t1} that any lsc strongly adequate mapping on a Banach
space $X$ is essentially firmly subdifferentiable and that the converse holds
for $X$ reflexive. In fact, Proposition \ref{cor3} says that the concept of
essentially firmly subdifferentiable mapping is intermediate between the
concept of essentially strongly convex function recently introduced in
\cite{Stromberg:11} and the concept of totally convex mapping (\cite{R5},
\cite{R6}, \cite{R7}, \cite{R14}, ...). In the reflexive case, the class of
lsc strongly adequate functions coincides with the one of essentially strongly
convex functions (Proposition \ref{prop7}). In the finite dimensional case,
the two classes above coincide with the essentially strictly convex functions
in the sense of \cite{R15} (Proposition \ref{cor3}). We provide an example of
an essentially strictly convex function on $\mathbb{R}^{2}$ with convex
subdifferential domain which is not totally convex (Example \ref{ex1}). A case
of essentially strictly convex function on $\mathbb{R}^{2}$ which is not
totally convex on the domain of its subdifferential is given in Example
\ref{ex2}.

An important tool we use is the natural notion of essential Fr\'{e}chet
differentiability introduced in \cite{Stromberg:11}, which strengthens the
concept of essential smoothness of \cite{R3}. In this way, a dual
characterization of lsc strongly adequate functions on general Banach spaces
is given in Proposition \ref{prop2} in terms of essentially Fr\'{e}chet
differentiability of the Legendre--Fenchel conjugate functions.

Section 3 is devoted to relative projections introduced in \cite{R7} which are
natural generalizations of the Bregman and Alber's projections (\cite{R1},
\cite{R4}, ...). In this context, we define the $f$-strongly Tchebychev sets
with respect to a lsc function $f$, non necessarily convex (see also
\cite{R3bis}), and study the convexity of this kind of sets (Theorem
\ref{t2}). As an application, we get a farthest point like result (Corollary
\ref{cor4}, Remark \ref{rem4}). In Section 4, we characterize the so-called
E-spaces (\cite{R8}, \cite{R9}, \cite{R11}, ...) by using our notion of
strongly adequate function (Proposition \ref{prop4}) and the firm
sub\-differentiability of the square of the norm (Proposition \ref{prop5}).
Finally, Proposition \ref{prop6} gives a variational characterization of the
closed convex sets in an E-space including the fact that a closed set in a
Hilbert space is convex if and only if it is strongly Tchebychev.

\section{Notation and preliminaries}

In the sequel $X$ is a Banach space whose topological dual and bidual we
denote by $X^{\ast}$ and $X^{\ast\ast}$, respectively; the dual norm on
$X^{\ast}$ is denoted by $\left\Vert \cdot\right\Vert _{\ast}$, and $x^{\ast
}(x)$ for $x\in X$ and $x^{\ast}\in X^{\ast}$ is denoted by $\left\langle
x,x^{\ast}\right\rangle $. We set $F(X)$ for the class of extended real-valued
functions $J:X\rightarrow\mathbb{R}\cup\{+\infty\}$ which are finite somewhere
(i.e.\ $\operatorname*{dom}J:=\{x\in X\mid J(x)<\infty\}\neq\emptyset$). As
usual $\Gamma(X)$ denotes the set of the lower semi\-continuous (lsc) convex
proper functions on $X$, and $J^{\ast}$ denotes the Legendre--Fenchel
conjugate of $J\in F(X):$
\[
J^{\ast}(x^{\ast}):=\sup_{x\in X}\left(  \left\langle x,x^{\ast}\right\rangle
-J(x)\right)  ,\quad x^{\ast}\in X^{\ast}.
\]

The subdifferential of $J\in F(X)$ ($J$ not necessarily convex) at a point
$x\in X$ is the set
\[
\partial J(x):=\{x^{\ast}\in X^{\ast}\mid J(u)\geq J(x)+\left\langle
u-x,x^{\ast}\right\rangle \ \forall u\in X\};
\]
clearly, $\partial J(x)=\emptyset$ if $J(x)\notin\mathbb{R}$. One has
\begin{equation}
x^{\ast}\in\partial J(x)\iff J^{\ast}(x^{\ast})+J(x)=\left\langle x,x^{\ast
}\right\rangle . \label{r-sub}%
\end{equation}
Consider the inverse multi\-map
\[
MJ:=(\partial J)^{-1}:X^{\ast}\rightrightarrows X;
\]
of course, one has%
\[
MJ(x^{\ast})=\arg\min(J-\left\langle \cdot,x^{\ast}\right\rangle ),\quad
x^{\ast}\in X^{\ast}.
\]
Taking into account (\ref{r-sub}) and the Fenchel--Young inequality one easily
observes that
\[
MJ(x^{\ast})\subset\partial J^{\ast}(x^{\ast})=\left\{  x^{\ast\ast}\in
X^{\ast\ast}\mid J^{\ast}(u^{\ast})\geq J^{\ast}(x^{\ast})+\left\langle
u^{\ast}-x^{\ast},x^{\ast\ast}\right\rangle \ \forall u^{\ast}\in X^{\ast
}\right\}  ,
\]
where $x\in X$ is identified with $\varphi\in X^{\ast\ast}$ defined by
$\varphi(x^{\ast}):=\left\langle x,x^{\ast}\right\rangle $. Since $J^{\ast}$
is sub\-differentiable on $\operatorname*{int}(\operatorname*{dom}J^{\ast})$
we have
\begin{equation}
\operatorname*{dom}MJ\cup\operatorname*{int}(\operatorname*{dom}J^{\ast
})\subset\operatorname*{dom}\partial J^{\ast}\subset\operatorname*{dom}%
J^{\ast}. \label{r1}%
\end{equation}

A mapping $J\in F(X)$ is said to be adequate (\cite{R16}) if
\[
\left\{
\begin{array}
[c]{l}%
\operatorname*{dom}MJ=\operatorname*{dom}(\partial J^{\ast})\text{ is a
nonempty open set, and}\\
MJ\text{ is single-valued on its domain.}%
\end{array}
\right.
\]
It has been proved in \cite[Th.\ 1]{R16} that for $X$ reflexive and $J\in
F(X)$, $J$ weakly lsc, $J$ is adequate iff $J$ is essentially strictly convex
in the sense of \cite{R3}.

\section{Strongly adequate functions}

Given an adequate function $J\in F(X)$ and $x^{\ast}\in\operatorname*{dom}MJ$,
the map $J-x^{\ast}$ attains a single minimum point over $X$. We are specially
interested in the case when this minimum is a strong minimum, that means every
minimizing sequence norm-converges to this minimum. To this end we recall
below an important result (see \cite[Cor.\ 6]{R2} and \cite[Prop.\ 4]{R12}):

\begin{lemma}
\label{lem1}Let $J\in F(X)$ be lsc and $x^{\ast}\in\operatorname*{int}%
(\operatorname*{dom}J^{\ast})$. Then $J-x^{\ast}$ attains a strong minimum
over $X$ iff $J^{\ast}$ is Fr\'{e}chet differentiable at $x^{\ast}.$
\end{lemma}

Such a situation occurs for instance in the following case:

\begin{proposition}
\label{prop1}Assume that $X$ has the Radon--Nikodym property and let $J\in
F(X)$ be lsc. Then the set of $x^{\ast}\in\operatorname*{int}%
(\operatorname*{dom}J^{\ast})$ such that $J-x^{\ast}$ attains a strong minimum
over $X$ is a dense $G_{\delta}$ in $\operatorname*{int}(\operatorname*{dom}%
J^{\ast}).$
\end{proposition}

Proof. By \cite[Th.\ 3.5.8]{R10}, $J^{\ast}$ is Fr\'{e}chet differentiable on
a dense $G_{\delta}$ subset $S$ of $\operatorname*{int}(\operatorname*{dom}%
J^{\ast})$. By Lemma \ref{lem1} $J-x^{\ast}$ attains a strong minimum over $X$
for every $x^{\ast}\in S$. \hfill$\square$

\medskip

In the light of the previous considerations we introduce now the notion of a
strongly adequate function.

\begin{definition}
\label{def1} A mapping $J\in F(X)$ is said to be \emph{strongly adequate} if
$\operatorname*{dom}\partial J^{\ast}$ is a nonempty open set and $J-x^{\ast}$
attains a strong minimum over $X$ for every $x^{\ast}\in\operatorname*{dom}%
\partial J^{\ast}.$
\end{definition}

According to (\ref{r1}) any strongly adequate $J\in F(X)$ satisfies%
\begin{equation}
\operatorname*{dom}MJ=\operatorname*{int}(\operatorname*{dom}J^{\ast
})=\operatorname*{dom}\partial J^{\ast}\neq\emptyset, \label{r2}%
\end{equation}
and, of course, any strongly adequate function is adequate.

An important example of strongly adequate function is furnished by the lsc
mappings $J\in F(X)$ whose conjugate $J^{\ast}$ is Fr\'{e}chet differentiable
on $X^{\ast}$. In order to go further in our investigation, let us quote the
following concept (see \cite[Def.\ 2]{Stromberg:11}).

\begin{definition}
\label{def2}Given a Banach space $Y$, we say that the function $G\in\Gamma(Y)$
is \emph{essentially Fr\'{e}chet differentiable} if $G$ is Fr\'{e}chet
differentiable at each point of $\operatorname*{int}(\operatorname*{dom}%
\partial G)\neq\emptyset$ and $\left\Vert \nabla G(x_{n})\right\Vert _{\ast
}\rightarrow\infty$ whenever $(x_{n})$ is a sequence in $\operatorname*{int}%
(\operatorname*{dom}G)$ converging to some boundary point of
$\operatorname*{dom}G.$
\end{definition}

\begin{proposition}
\label{p-stromberg}A mapping $G\in\Gamma(Y)$ is essentially Fr\'{e}chet
differentiable iff $G$ is Fr\'{e}chet differentiable at each point of
$\operatorname*{dom}\partial G.$
\end{proposition}

Proof. It is similar to that of the equivalence of (i) and (v) in
\cite[Th.\ 5.6]{R3}. We give the proof for reader's convenience.

Sufficiency: here $\operatorname*{dom}\partial G$ is open, and so
$\operatorname*{dom}\partial G=\operatorname*{int}(\operatorname*{dom}G)$.
Hence $G$ is Fr\'{e}chet differentiable on $\operatorname*{int}%
(\operatorname*{dom}G)$. Let $(x_{n})\subset\operatorname*{int}%
(\operatorname*{dom}G)$ be convergent to $x\in\operatorname*{bd}%
(\operatorname*{dom}G)$, and assume that $\left\Vert \nabla G(x_{n}%
)\right\Vert _{\ast}\not \rightarrow \infty$. Passing to a subsequence if
necessary, we may (and do) assume that $\big(
\left\Vert \nabla G(x_{n})\right\Vert _{\ast}\big)  $ is bounded. Therefore,
$\big( \left\Vert \nabla G(x_{n})\right\Vert _{\ast}\big)
$ has a subnet $\big( \left\Vert \nabla G(x_{\varphi(i)})\right\Vert _{\ast
}\big)  _{i\in I}$ converging weakly-star to $x^{\ast}\in Y^{\ast}$. By
\cite[Th.\ 2.4.2(ix)]{R17} we obtain that $(x,x^{\ast})\in\partial G$, and so
we get the contradiction $x\in\operatorname*{dom}\partial
G=\operatorname*{int} (\operatorname*{dom}G)$. Therefore, $G$ is essentially
Fr\'{e}chet differentiable.

Necessity: since $G$ is essentially Fr\'{e}chet differentiable on
$\operatorname*{int}(\operatorname*{dom}G)\neq\emptyset$, we have to show that
$\operatorname*{dom}\partial G=\operatorname*{int}(\operatorname*{dom}G)$.
Assume that there exists $x\in\operatorname*{dom}\partial G\setminus
\operatorname*{int}(\operatorname*{dom}G)$. Fix $\overline{x}\in
\operatorname*{int}(\operatorname*{dom}G);$ clearly, $]x,\overline{x}%
]\subset\operatorname*{int}(\operatorname*{dom}G)$. Using \cite[Lem.\ 4.4]%
{R3}, we have that $\nabla G(]x,\overline{x}])$ is bounded. Taking
$x_{n}:=(1-n^{-1})x+n^{-1}\overline{x}\in{}]x,\overline{x}]$, we have that
$x_{n}\rightarrow x$ and $\left(  \nabla G(x_{n})\right)  $ is bounded. This
contradiction proves that $\operatorname*{dom}\partial G\subset
\operatorname*{int}(\operatorname*{dom}G)$. \hfill$\square$

\medskip

An essentially Fr\'{e}chet differentiable function $G\in\Gamma(Y)$ satisfies
$\operatorname*{dom}\partial G=\operatorname*{int}(\operatorname*{dom}%
G)\neq\emptyset$ and $\partial G$ is both single-valued and locally bounded on
its domain (see e.g.\ \cite[Cor.\ 2.4.13]{R17}). Consequently, any essentially
Fr\'{e}chet differentiable function $G\in\Gamma(Y)$ is essentially smooth in
the sense of \cite[Def.\ 5.2]{R3}.

It is worthwhile noting that if $Y$ is finite dimensional the two notions
above coincide with the usual one introduced in \cite[Section 2.6]{R15}.

\begin{remark}
\label{rem1}If $G\in\Gamma(Y)$ is finite-valued, then $G$ is essentially
Fr\'{e}chet differentiable iff $G$ is Fr\'{e}chet differentiable at each point
of $Y.$
\end{remark}

We now provide a dual characterization for a strongly adequate function.

\begin{proposition}
\label{prop2}A lsc mapping $J\in F(X)$ is strongly adequate iff its conjugate
$J^{\ast}$ is essentially Fr\'{e}chet differentiable.
\end{proposition}

Proof. In both cases $\operatorname*{dom}\partial J^{\ast}$ is open, nonempty
and, according to (\ref{r1}) and (\ref{r2}), coincides with
$\operatorname*{int}(\operatorname*{dom}J^{\ast})$. It then suffices to apply
Lemma \ref{lem1}. \hfill$\square$

\medskip

As in \cite{R1bis} (see also \cite[p.\ 188]{R17} and \cite{Stromberg:11}), let
us introduce the set
\[
\Gamma_{0}:=\left\{  \psi:\mathbb{R}_{+}\rightarrow\lbrack0,\infty]\mid
\psi\text{ lsc convex, }\psi(t)=0\iff t=0\right\}  .
\]
Any $\psi\in\Gamma_{0}$ is a forcing function in the sense of \cite[p.\ 6]%
{R8}:%
\[
\forall(t_{n})\subset\mathbb{R}_{+}:\psi(t_{n})\rightarrow0\Rightarrow
t_{n}\rightarrow0.
\]
Also, any $\psi\in\Gamma_{0}$ satisfies $\lim_{t\rightarrow\infty}%
\psi(t)=\infty.$

The following concept has been introduced in \cite[Def.\ 2]{Stromberg:11}.

\begin{definition}
\label{def-sc}A mapping $H\in\Gamma(X)$ is said to be \emph{essentially
strongly convex} if it is essentially strictly convex in the sense of
\cite[Def.\ 5.2]{R3} and if for every $x\in\operatorname*{dom}\partial H$
there exist $x^{\ast}\in\partial H(x)$ and $\psi\in\Gamma_{0}$ such that:
\begin{equation}
H(u)\geq H(x)+\left\langle u-x,x^{\ast}\right\rangle +\psi(\left\Vert
u-x\right\Vert )\quad\forall u\in X. \label{r3}%
\end{equation}

\end{definition}

By \cite[Th.\ 3]{Stromberg:11} we know that for any $J\in F(X)$, one has:%
\[
J^{\ast}\text{ essentially Fr\'{e}chet differentiable }\Rightarrow\text{
}J\text{ essentially strongly convex.}%
\]

>From our Proposition \ref{prop2} above we thus have:

\begin{corollary}
\label{cor1}Any lsc strongly adequate function $J\in F(X)$ is essentially
strongly convex.
\end{corollary}

In fact, more can be said. To this end, let us introduce the following notion
(which appears in \cite[Prop.\ 2]{Stromberg:11} in the framework of
essentially strictly convex functions on reflexive Banach spaces):

\begin{definition}
\label{def3}A convex mapping $H\in F(X)$ is said to be \emph{firmly
sub\-differentiable at} $x\in\operatorname*{dom}\partial H$ if for any
$x^{\ast}\in\partial H(x)$ there exists $\psi\in\Gamma_{0}$ such that
$(\ref{r3})$ holds. If $H$ is firmly sub\-differentiable at each point of
$\operatorname*{dom}\partial H$ we will say that $H$ is \emph{essentially
firmly sub\-differentiable}.
\end{definition}

In order to illustrate Definition \ref{def3} let us recall that a convex
mapping $H\in F(X)$ is said to be totally convex at a point $x\in
\operatorname*{dom}H$ if, denoting by $H^{\prime}(x,d)$ the right hand side
derivative of $H$ at $x$ in the direction $d$, one has (\cite{R5})%
\[
\inf\left\{  H(u)-H(x)-H^{\prime}(x,u-x)\mid u\in\operatorname*{dom}%
H,\ \left\Vert u-x\right\Vert =t\right\}  >0\quad\forall t>0.
\]
Given $x\in\operatorname*{dom}\partial H$ we know (\cite[Lem.\ 3.3]{R6}) that
$H$ is totally convex at $x$ iff there exists $\xi\in\Gamma_{0}$ such that
($\operatorname*{int}(\operatorname*{dom}\xi)\neq\emptyset$ and)%
\[
H(u)\geq H(x)+H^{\prime}(x,u-x)+\xi(\left\Vert u-x\right\Vert )\quad\forall
u\in X.
\]
Since for any $x^{\ast}\in\partial H(x)$ and any direction $d$ one has
$\left\langle d,x^{\ast}\right\rangle \leq H^{\prime}(x,d)$, it holds that if
$H$ is totally convex at $x\in\operatorname*{dom}\partial H$ then $H$ is
firmly sub\-differentiable at $x.$

It follows from \cite[Prop.\ 3.5]{R6} that for the proper convex function $H$
which is continuous at $x\in\operatorname*{dom}H$ we have that $H$ is totally
convex at $x$ iff $H$ is uniformly firmly sub\-differentiable at $x$ (that is
the same $\psi$ is valid for all $x^{\ast}\in\partial H(x)$).

\begin{proposition}
\label{prop-z}Let $X$ be finite dimensional, $H\in F(X)$, $H$ convex, and
$\overline{x}\in\operatorname*{rint}(\operatorname*{dom}H)$. Then $H$ is
totally convex at $\overline{x}$ iff $H$ is firmly sub\-differentiable at
$\overline{x}.$
\end{proposition}

Proof. Clearly $\partial H(\overline{x})\neq\emptyset$. We may (and do) assume
that $H$ is lsc. The implication $\Rightarrow$ was observed above. Replacing
$H$ by $H(\overline{x}+\cdot)-H(\overline{x})$ we may assume that
$\overline{x}=0$ and $H(0)=0$. Moreover, taking $X_{0}=\operatorname*{aff}%
(\operatorname*{dom}H)=\operatorname*{lin}(\operatorname*{dom}H)$, we have
that $\operatorname*{dom}H^{\prime}(0,\cdot)=X_{0}$. It follows that
$H^{\prime}(0,\cdot)|_{X_{0}}$ is continuous on $X_{0}$. Assume that $H$ is
firmly sub\-differentiable at $\overline{x}=0$ but $H$ is not totally convex
at $\overline{x}$. Then there exists $\overline{t}>0$ such that
\[
\inf\left\{  H(x)-H^{\prime}(0,x)\mid x\in\operatorname*{dom}H,\ \left\Vert
x\right\Vert =\overline{t}\right\}  =\inf\left\{  H(x)-H^{\prime}(0,x)\mid
x\in X_{0},\ \left\Vert x\right\Vert =\overline{t}\right\}  =0.
\]
Because $H|_{X_{0}}-H^{\prime}(0,\cdot)|_{X_{0}}$ is lsc on $X_{0}$ and
$A:=\left\{  x\in X_{0}\mid\left\Vert x\right\Vert =\overline{t}\right\}  $ is
compact, there exists $\overline{u}\in A$ such that $H(\overline{u}%
)-H^{\prime}(0,\overline{u})=0$. But $H^{\prime}(0,u)=\max\{\left\langle
u,x^{\ast}\right\rangle \mid x^{\ast}\in\partial H(0)\}$ for every $u\in
X_{0}$, and so there exists $\overline{x}^{\ast}\in\partial H(0)$ with
$H^{\prime}(0,\overline{u})=\left\langle \overline{u},\overline{x}^{\ast
}\right\rangle $. Therefore, $\inf\left\{  H(x)-\left\langle x,\overline
{x}^{\ast}\right\rangle \mid\left\Vert x\right\Vert =\overline{t}\right\}
=0$, contradicting the fact that $H$ is firmly sub\-differentiable at $0$.
\hfill$\square$

\medskip

In \cite{R14} it was posed the problem if, in finite dimensions, the converse
of \cite[Prop.\ 2.1]{R14} is true, that is if an essentially strictly convex
function $H\in\Gamma(X)$ with $\operatorname*{dom}\partial H$ convex is
totally convex (in the case $\dim X<\infty)$. We provide an example of an
essentially strictly convex function $H\in\Gamma(\mathbb{R}^{2})$ with
$\operatorname*{dom}\partial H$ convex which is not totally convex.

\begin{example}
\label{ex1}Let $H:\mathbb{R}^{2}\rightarrow\overline{\mathbb{R}}$ be defined
by $H(x,y):=-\sqrt[4]{(1-x^{2})(1-y^{2})}$ for $(x,y)\in\lbrack-1,1]\times
\lbrack-1,1]$, $H(x,y):=\infty$ otherwise. Observe that
$H|_{\operatorname*{dom}H}$ is continuous and%
\begin{gather*}
\frac{\partial^{2}H(x,y)}{\partial x^{2}}=\frac{x^{2}+2}{4}\sqrt[4]%
{\frac{1-y^{2}}{\left(  1-x^{2}\right)  ^{7}}}>0,\\
\frac{\partial^{2}H(x,y)}{\partial x^{2}}\frac{\partial^{2}H(x,y)}{\partial
y^{2}}-\left(  \frac{\partial^{2}H(x,y)}{\partial x\partial y}\right)
^{2}=\frac{x^{2}+y^{2}+2}{8[\left(  1-x^{2}\right)  \left(  1-y^{2}\right)
]^{3/2}}>0
\end{gather*}
on $(-1,1)\times(-1,1)$. It follows that $H$ is convex, Fr\'{e}chet
differentiable and strictly convex on $\operatorname*{dom}\partial
H=(-1,1)\times(-1,1)$, and $\operatorname*{dom}H^{\ast}=\mathbb{R}^{2}$. It
follows that $H$ is essentially strictly convex. Since $H$ is not strictly
convex ($H(x,1)=0$ for every $x\in\lbrack-1,1]$), we have that $H$ is not
totally convex.
\end{example}

Let us provide some properties of essentially firmly sub\-differentiable
mappings in terms of well-posedness and coercivity.

\begin{proposition}
\label{prop5b}Let $H$ be firmly sub\-differentiable at a point of
$\operatorname{argmin}H$ (supposed to be nonempty). Then $H$ is coercive and
attains a strong minimum over $X.$
\end{proposition}

Proof. Let $x\in\operatorname{argmin}H$ a point where $H$ is firmly
sub\-differentiable. One has $0\in\partial H(x)$ and there exists $\psi
\in\Gamma_{0}$ such that
\begin{equation}
H(u)\geq H(x)+\psi(\left\Vert u-x\right\Vert ),\quad\forall u\in X. \label{r4}%
\end{equation}
Let $(x_{n})$ be a minimizing sequence of $H$. By (\ref{r4}) we have that
$\psi\left(  \left\Vert x_{n}-x\right\Vert \right)  \rightarrow0$ and so
$\left\Vert x_{n}-x\right\Vert \rightarrow0$. One has also $\lim
_{t\rightarrow\infty}\psi(t)=\infty$ and by (\ref{r4}) $\lim_{\left\Vert
u\right\Vert \rightarrow\infty}H(u)=\infty$. \hfill$\square$

\medskip

We now appeal to a dual interpretation of relation (\ref{r3}):

\begin{lemma}
[\cite{R17}, Cor.\ 3.4.4]\label{lem2} Let $H\in\Gamma(X)$ and $(x,x^{\ast}%
)\in\partial H$. The statements below are equivalent:

\emph{(i)} $\exists\psi\in\Gamma_{0}$ such that $(\ref{r3})$ holds,

\emph{(ii)} $H^{\ast}$ is Fr\'{e}chet differentiable at $x^{\ast}.$
\end{lemma}

We are now in a position to characterize the strongly adequate mappings among
the lsc ones:

\begin{theorem}
\label{t1}Let $J\in F(X)$ be lsc. If $J$ is strongly adequate then $J$ is
essentially firmly sub\-differentiable. The converse holds for $X$ reflexive.
\end{theorem}

Proof. Assume $J$ is strongly adequate. By Corollary \ref{cor1} we have that
$J\in\Gamma(X)$. Let $(x,x^{\ast})\in\partial J$. Thus $x^{\ast}%
\in\operatorname*{dom}\partial J^{\ast}$, and Proposition \ref{prop2} says
that $J^{\ast}$ is Fr\'{e}chet differentiable at $x^{\ast}$. By Lemma
\ref{lem2} there exists $\psi\in\Gamma_{0}$ such that (\ref{r3}) holds,
meaning that $J$ is firmly sub\-differentiable at any point $x\in
\operatorname*{dom}\partial J$, thus essentially firmly sub\-differentiable.

Assume $J$ is essentially firmly sub\-differentiable and $X$ is reflexive. Let
$x^{\ast}\in\operatorname*{dom}\partial J^{\ast}$. Since $X$ is reflexive
there is $x\in X$ such that $(x,x^{\ast})\in\partial J$. Since $J$ is firmly
sub\-differentiable at $x$, Lemma \ref{lem2} says that $J^{\ast}$ is
Fr\'{e}chet differentiable at $x^{\ast}$. Consequently, $J^{\ast}$ is
essentially Fr\'{e}chet differentiable and, by Proposition \ref{prop2}, $J$ is
strongly adequate. \hfill$\square$

\medskip

Let us establish some links among some of the convexity notions quoted above.

\begin{proposition}
\label{cor3}Let $H\in\Gamma(X)$, and consider the following statements:

\emph{(i)} $H$ is totally convex at each point of $\operatorname*{dom}\partial
H,$

\emph{(ii)} $H$ is essentially firmly sub\-differentiable,

\emph{(iii)} $H$ is essentially strongly convex,

\emph{(iv)} $H$ is essentially strictly convex.

Then \emph{(i)} $\Rightarrow$ \emph{(ii)} $\Rightarrow$ \emph{(iii)}
$\Rightarrow$ \emph{(iv)}. Moreover, if $X$ is finite dimensional then
\emph{(iv)} $\Rightarrow$ \emph{(ii)}.
\end{proposition}

Proof. (i) $\Rightarrow$ (ii) was observed after Definition \ref{def3}.

(ii) $\Rightarrow$ (iv) Let $x_{0},x_{1}\in X$ be such that $[x_{0}%
,x_{1}]\subset\operatorname*{dom}\partial H$. Assume that $H$ is not strictly
convex on $[x_{0},x_{1}];$ then $x_{0}\neq x_{1}$ and there exists $\lambda
\in(0,1)$ such that $H(x_{\lambda})=(1-\lambda)H(x_{0})+\lambda H(x_{1})$,
where $x_{\lambda}:=(1-\lambda)x_{0}+\lambda x_{1}$. Take $x^{\ast}\in\partial
H(x_{\lambda})$ and $\psi\in\Gamma_{0}$ for which (\ref{r3}) holds. Then%
\[
H(x_{i})\geq H(x_{\lambda})+\left\langle x_{i}-x_{\lambda},x^{\ast
}\right\rangle +\psi(\left\Vert x_{i}-x_{\lambda}\right\Vert )\quad
(i\in\{0,1\}).
\]
Multiplying both terms of this inequality by $1-\lambda>0$ for $i=0$ and by
$\lambda>0$ for $i=1$, then adding side by side we get%
\[
H(x_{\lambda})=(1-\lambda)H(x_{0})+\lambda H(x_{1})\geq H(x_{\lambda
})+(1-\lambda)\psi(\left\Vert x_{0}-x_{\lambda}\right\Vert )+\lambda
\psi(\left\Vert x_{1}-x_{\lambda}\right\Vert ).
\]
It follows that $\psi(\left\Vert x_{0}-x_{\lambda}\right\Vert )=\psi
(\left\Vert x_{1}-x_{\lambda}\right\Vert )=0$, whence the contradiction
$x_{0}=x_{1}$ $(=x_{\lambda}).$

Let now $x^{\ast}\in\operatorname*{dom}(\partial H)^{-1}$. Then there exists
$x\in\operatorname*{dom}\partial H$ such that $x^{\ast}\in\partial H(x)$. By
(ii), there exists $\psi\in\Gamma_{0}$ such that (\ref{r3}) holds. Using Lemma
\ref{lem2} we obtain that $H^{\ast}$ is Fr\'{e}chet differentiable at
$x^{\ast}$, and so $x^{\ast}\in\operatorname*{int}(\operatorname*{dom}H^{\ast
})$; hence $\partial H^{\ast}$ is bounded on a neighborhood $V$ of $x^{\ast}$.
Because $(\partial H)^{-1}(u^{\ast})\subset\partial H^{\ast}(u^{\ast})$ for
every $u^{\ast}\in X^{\ast}$, we obtain that $(\partial H)^{-1}$ is bounded on
$V.$

(ii) $\Rightarrow$ (iii) follows immediately from (ii) $\Rightarrow$ (iv) and
the fact that $\partial H(x)\neq\emptyset$ for every $x\in\operatorname*{dom}%
\partial H.$

(iv) $\Rightarrow$ (ii) Let $\dim X<\infty$. By \cite[Th.\ 2.6.3]{R15} (or
\cite[Th.\ 5.4]{R3}) we know that $H^{\ast}$ is essentially smooth. Since $X$
is finite dimensional this amounts to say that $H^{\ast}$ is essentially
Fr\'{e}chet differentiable, and, by Corollary \ref{cor2}, that $H$ is
essentially firmly sub\-differentiable. \hfill$\square$

\medskip

We have seen in Example \ref{ex1} that, for $\dim X<\infty$, there exist
essentially strictly convex functions which are not totally convex. A natural
question is if (iv) $\Rightarrow$ (i) in Corollary \ref{cor3} for $\dim
X<\infty$. The answer is negative, as shown in the following example.

\begin{example}
\label{ex2}Let $H:\mathbb{R}^{2}\rightarrow\overline{\mathbb{R}}$ be defined
by $H(x,y):=-\sqrt{(1-x^{2})(1-y^{2})}$ for $(x,y)\in\lbrack-1,1]\times
\lbrack-1,1]$, $H(x,y):=\infty$ otherwise. Then $\operatorname*{dom}\partial
H=[(-1,1)\times(-1,1)]\cup\lbrack\{-1,1\}\times\{-1,1\}]$, $H$ is strictly
convex on every segment included in $\operatorname*{dom}\partial H$ and
$\operatorname*{dom}H^{\ast}=\mathbb{R}^{2}$. It follows that $H$ is
essentially strictly convex. Since $H$ is not strictly convex ($H(x,1)=0$ for
every $x\in\lbrack-1,1]$), we have that $H$ is not totally convex at each
point of $\operatorname*{dom}\partial H.$
\end{example}

\begin{remark}
\label{rem3}Proposition \ref{cor3} provides in particular a significant
improvement of \cite[Prop.\ 2.1]{R14}.
\end{remark}

By juxtaposition of Theorem \ref{t1} and Proposition \ref{prop2} we obtain:

\begin{corollary}
\label{cor2}Let $J\in F(X)$ be lsc. If $J^{\ast}$ is essentially Fr\'{e}chet
differentiable, then $J$ is essentially firmly sub\-differentiable. The
converse holds for $X$ reflexive.
\end{corollary}

\begin{remark}
\label{rem-z}According to Proposition \ref{cor3}, the first part of Corollary
\ref{cor2} improves \cite[Th.\ 3]{Stromberg:11}.
\end{remark}

\begin{remark}
According to Remark \ref{rem1} and Corollary \ref{cor2}, a co\-finite mapping
$H\in\Gamma(X)$ (that is $H^{\ast}$ is finite-valued) with $X$ reflexive is
essentially firmly sub\-differentiable iff $H^{\ast}$ is essentially
Fr\'{e}chet differentiable on $X^{\ast}$. For instance, the square of the norm
of a reflexive Banach space is essentially firmly sub\-differentiable iff the
square of the dual norm is Fr\'{e}chet differentiable on $X^{\ast}$ (see
Section 4).
\end{remark}

We end this section by a more complete result concerning the reflexive case.
It includes \cite[Th.\ 4]{Stromberg:11} and a part of \cite[Prop.\ 2]%
{Stromberg:11}.

\begin{proposition}
\label{prop7} Let $X$ be a reflexive Banach space and $H\in F(X)$. The
following assertions are equivalent:

\emph{(i)} $H$ is strongly adequate,

\emph{(ii)} $H^{\ast}$ is essentially Fr\'{e}chet differentiable,

\emph{(iii)} $H$ is essentially firmly sub\-differentiable,

\emph{(iv)} $H$ is essentially strongly convex.
\end{proposition}

Proof. (i)\ $\Leftrightarrow$ (ii) is established in Theorem \ref{t1}, and
(ii)\ $\Leftrightarrow$ (iii) in Corollary \ref{cor2}. By Proposition
\ref{cor3} one has (iii) $\Rightarrow$ (iv). The equivalence
(ii)\ $\Leftrightarrow$ (iv) is \cite[Th.\ 4]{Stromberg:11}.

We give a proof of (iv) $\Rightarrow$ (ii) for reader's convenience. By
Definition \ref{def-sc} we have that $H$ is essentially strictly convex, and
so, by \cite[Th.\ 5.4]{R3}, $H^{\ast}$ is essentially smooth. By
\cite[Th.\ 5.6]{R3} we obtain that $\operatorname*{dom}\partial H^{\ast
}=\operatorname*{int}(\operatorname*{dom}H^{\ast})\neq\emptyset$ and $H^{\ast
}$ is G\^{a}teaux differentiable on $\operatorname*{int}(\operatorname*{dom}%
H^{\ast})$. Let $\overline{x}^{\ast}\in\operatorname*{int}(\operatorname*{dom}%
H^{\ast})$ and let us show that $H^{\ast}$ is Fr\'{e}chet differentiable at
$\overline{x}^{\ast}$. For this we apply \cite[Th.\ 3.3.2]{R17}. Set
$\overline{x}:=\nabla H^{\ast}(\overline{x}^{\ast})\in\operatorname*{dom}%
\partial H$ and take $\left(  (x_{n},x_{n}^{\ast})\right)  _{n\geq1}%
\subset\partial H$ with $x_{n}^{\ast}\rightarrow\overline{x}^{\ast}$. By
\cite[Th.\ 3.3.2]{R17} applied for the G\^{a}teaux bornology we have that
$x_{n}\rightarrow^{w}\overline{x}$. Since $H(\overline{x})\geq H(x_{n}%
)+\left\langle \overline{x}-x_{n},x_{n}^{\ast}\right\rangle $ for every
$n\geq1$, we obtain that $H(\overline{x})\geq\limsup H(x_{n})$, and so
$H(x_{n})\rightarrow H(\overline{x})$ because $H$ is weakly lsc. Because $H$
is essentially strongly convex, there exist $\widetilde{x}^{\ast}\in\partial
H(\overline{x})$ and $\psi\in\Gamma_{0}$ such that $H(x)\geq H(\overline
{x})+\left\langle x-\overline{x},\widetilde{x}^{\ast}\right\rangle
+\psi(\left\Vert x-\overline{x}\right\Vert )$ for every $x\in X$, and so
$H(x_{n})\geq H(\overline{x})+\left\langle x_{n}-\overline{x},\widetilde
{x}^{\ast}\right\rangle +\psi(\left\Vert x_{n}-\overline{x}\right\Vert )$ for
every $n\geq1$. Taking the $\limsup$ in both terms we get $\limsup
\psi(\left\Vert x_{n}-\overline{x}\right\Vert )\leq0$, and so $\left\Vert
x_{n}-\overline{x}\right\Vert \rightarrow0$. Applying now \cite[Th.\ 3.3.2]%
{R17} for the Fr\'{e}chet bornology we obtain that $H^{\ast}$ is Fr\'{e}chet
differentiable at $\overline{x}^{\ast}$. \hfill$\square$

\section{Relative projections on closed sets}

Given $f\in F(X)$, $S$ a closed subset of the Banach space $X$ such that
\[
S\cap\operatorname*{dom}f\neq\emptyset,
\]
and $x^{\ast}\in X^{\ast}$, let us consider the problem:
\[
P_{S}(f,x^{\ast}):\quad\min\left(  f(x)-\left\langle x,x^{\ast}\right\rangle
\right)  \quad\text{for }x\in S.
\]

Such problems have been studied in \cite{R7} under the name relative
projection (of $x^{\ast}$ on $S$ modulo $f$). They are natural generalizations
of the Bregman projections and generalized projections defined and studied by
Alber (\cite{R1}, \cite{R4}, ...). In \cite{R7} the mapping $f$ is assumed to
be convex. We don't retain this assumption here, and just assume that $f\in
F(X)$. For instance, taking $f:=-\tfrac{1}{2}\left\Vert \cdot\right\Vert ^{2}$
on the Hilbert space $(X,\left\Vert \cdot\right\Vert )$ and $S\subset X$ a
bounded subset, the problem $P_{S}(f,x^{\ast})$ consists of finding the
farthest points of $S$ from $-x^{\ast}\in X^{\ast}=X$. Taking $f:=\tfrac{1}%
{2}\left\Vert \cdot\right\Vert ^{2}$ and $S\subset X$, still in the Hilbert
space setting, the problem $P_{S}(f,x^{\ast})$ consists in finding the best
approximation of $x^{\ast}\in X^{\ast}=X$ by elements of $S.$

\begin{definition}
\label{def4}We will say that $S$ is $f$\emph{-strongly Tchebychev} if for
every $x^{\ast}\in X^{\ast}$ the problem $P_{S}(f,x^{\ast})$ admits a strong
minimum; in other words if any minimizing sequence of $P_{S}(f,x^{\ast})$
norm-converges toward a (necessarily unique) solution of $P_{S}(f,x^{\ast}).$
\end{definition}

Denoting by $\iota_{S}$ the indicator function of $S$, it is clear that if $S$ is $f$-strongly Tchebychev, then
$f+\iota_{S}$ is strongly adequate; conversely, if $f+\iota_{S}$ is  strongly adequate and cofinite, then $S$ is $f$-strongly Tchebychev. We can state:

\begin{theorem}
\label{t2}Let $X$ be a Banach space, $f\in F(X)$, $f$ lsc, and $S\subset X$
closed satisfying $S\cap\operatorname*{dom}f\neq\emptyset$. If $S$ is
$f$-strongly Tchebychev then $f+\iota_{S}$ is essentially firmly
sub\-differentiable and $S\cap\operatorname*{dom}f$ is convex.

Conversely, if $X$ is reflexive, $S$ is convex, and $f$ is essentially firmly
sub\-differentiable, finite and continuous at a point of $S$, then
$f+\iota_{S}$ is strongly adequate; moreover, if $f+\iota_{S}$ is cofinite
then $S$ is $f$-strongly Tchebychev.
\end{theorem}

Proof. Assume $S$ is $f$-strongly Tchebychev. Then $J:=f+\iota_{S}$ is
strongly adequate and, by Theorem \ref{t1}, $J$ is essentially firmly
sub\-differentiable. In particular, $J$ is convex and so $\operatorname*{dom}%
J=S\cap\operatorname*{dom}f$ is convex.

Conversely, let us first notice that, by \cite[Prop.\ 10.d]{R13} or
\cite[Th.\ 2.8.7(iii)]{R17}, one has $\partial(f+\iota_{S})(x)=\partial
f(x)+\partial\iota_{S}(x)$ for all $x\in X$. We thus have $\operatorname*{dom}%
\partial J=S\cap\operatorname*{dom}\partial f$. Now for any $x\in
\operatorname*{dom}\partial J$, any $x^{\ast}\in\partial J(x)$, there exist
$u^{\ast}\in\partial f(x)$ and $v^{\ast}\in N(S,x)$ such that $x^{\ast
}=u^{\ast}+v^{\ast}$. Since $f$ is firmly subdifferentiable at $x$, there
exists $\psi\in\Gamma_{0}$ such that, for any $u\in X$, $f(u)\geq
f(x)+\left\langle u-x,u^{\ast}\right\rangle +\psi(\left\Vert u-x\right\Vert
)$, and thus $J(u)\geq J(x)+\left\langle u-x,x^{\ast}\right\rangle
+\psi(\left\Vert u-x\right\Vert )$, that means $J$ is firmly subdifferentiable
at each $x\in\operatorname*{dom}\partial J$. By the second part of Theorem
\ref{t1} we infer that $J$ is strongly adequate. When  $f+\iota_{S}$ is
cofinite this means that $S$ is $f$-strongly Tchebychev. \hfill$\square$

\begin{corollary}
\label{cor4}Let $X$ be a Banach space and $S$ a nonempty closed bounded subset
of $X$. The following statements are equivalent:

\emph{(i)} the mapping $\tfrac{1}{2}\left\Vert \cdot\right\Vert ^{2}+x^{\ast}$
attains a strong maximum over $S$ for any $x^{\ast}\in X^{\ast}$,

\emph{(ii)} $S$ is a singleton.
\end{corollary}

Proof. It is clear that (ii) $\Rightarrow$ (i). Conversely, (i) says that the
function $J:=-\tfrac{1}{2}\left\Vert \cdot\right\Vert ^{2}+\iota_{S}$ is
strongly adequate, hence convex by Corollary \ref{cor1}. It follows that $S$
is convex. Because $J$ is strongly adequate, there exists $x_{0}\in S$ such
that $J(x_{0})<J(x)$ for every $x\in S\setminus\{x_{0}\}$. Assume that
$S\neq\{x_{0}\}$, and take $x_{1}\in S\setminus\{x_{0}\}$. Then $x_{\lambda
}:=(1-\lambda)x_{0}+\lambda x_{1}\in S\setminus\{x_{0}\}$ and
\[
-\tfrac{1}{2}\left\Vert x_{\lambda}\right\Vert ^{2}=J(x_{\lambda}%
)\leq(1-\lambda)J(x_{0})+\lambda J(x_{1})=-(1-\lambda)\tfrac{1}{2}\left\Vert
x_{0}\right\Vert ^{2}-\lambda\tfrac{1}{2}\left\Vert x_{1}\right\Vert ^{2}%
\]
for every $\lambda\in\lbrack0,1]$. Since $\left\Vert x_{\lambda}\right\Vert
^{2}\leq(1-\lambda)\left\Vert x_{0}\right\Vert ^{2}+\lambda\left\Vert
x_{1}\right\Vert ^{2}$ for $\lambda\in\lbrack0,1]$ with strict inequality for
$\lambda\in(0,1)$ and $\left\Vert x_{0}\right\Vert \neq\left\Vert
x_{1}\right\Vert ,$ we obtain that $\left\Vert x_{0}\right\Vert =\left\Vert
x_{1}\right\Vert $, and so we get the contradiction $J(x_{0})=J(x_{1})$. Hence
$S$ is necessarily a singleton. \hfill$\square$

\begin{remark}
\label{rem4}In the Hilbert space setting, Corollary \ref{cor4} gives the
equivalence between the next two statements (involving farthest points to the
nonempty closed bounded set $S\subset X$):

\emph{(i)} for any $x^{\ast}\in X^{\ast}=X$, the mapping $x\longmapsto
\left\Vert x^{\ast}-x\right\Vert $ attains a strong maximum over $S,$

\emph{(ii)} $S$ is a singleton.
\end{remark}

\section{Variational characterizations of closed convex sets in E-spaces}

Let us recall that a Banach space $X$ is is said to be an E-space if $X$ is
rotund and every weakly closed set in $X$ is approximately compact. Such
spaces, introduced in \cite{R9}, admit several characterizations. For
instance, the theorem in \cite[p.\ 146]{R11} says that $X$ is an E-space iff
$X$ is reflexive, rotund, and any weakly convergent sequence within the unit
sphere of $X$ is norm convergent. Anderson's Theorem (see \cite[p.\ 149]{R11})
says that the Banach space $X$ is an E-space iff the square of the dual norm
$\left\Vert \cdot\right\Vert _{\ast}^{2}$ is Fr\'{e}chet differentiable on
$X^{\ast}$. In the light of the previous results we thus can obtain other
characterizations for the E-spaces.

\begin{proposition}
\label{prop4}For any Banach space $X$, the statements below are equivalent:

\emph{(i)} For any $x^{\ast}\in X^{\ast}$, the mapping $\tfrac{1}{2}\left\Vert
\cdot\right\Vert ^{2}-x^{\ast}$ attains a strong minimum over $X,$

\emph{(ii)} $X$ is an E-space.
\end{proposition}

Proof. Setting $J:=\tfrac{1}{2}\left\Vert \cdot\right\Vert ^{2}$, one has
$J^{\ast}=\tfrac{1}{2}\left\Vert \cdot\right\Vert _{\ast}^{2};$ the
equivalence between (i) and (ii) is obtained by Proposition \ref{prop2} and
Remark \ref{rem1} using Anderson's Theorem mentioned above. \hfill
\hfill$\square$

\medskip

In \cite[Th.\ 3.3]{R14}, the E-spaces are characterized among the reflexive
Banach space as the locally totally convex spaces (that are the reflexive
Banach spaces whose square of the norm is totally convex at any point). Below
we characterize the E-spaces in terms of the firm sub\-differentiability of
the square of the norm:

\begin{proposition}
\label{prop5}Given a Banach space $(X,\left\Vert \cdot\right\Vert )$, let us
consider the following statements:

\emph{(i)} $X$ is an E-space,

\emph{(ii)} $\left\Vert \cdot\right\Vert ^{2}$ is totally convex on $X,$

\emph{(iii)} $\left\Vert \cdot\right\Vert ^{2}$ is essentially firmly sub\-differentiable.

Then we have \emph{(i)} $\Rightarrow$ \emph{(ii)} $\Rightarrow$ \emph{(iii)},
and, if $X$ is reflexive, then \emph{(i)} $\Leftrightarrow$ \emph{(ii)}
$\Leftrightarrow$ \emph{(iii)}.
\end{proposition}

Proof. (i) $\Rightarrow$ (ii) is proved in \cite[Th.\ 3.3]{R14}, and (ii)
$\Rightarrow$ (iii) has been quoted in the comments after Definition
\ref{def3}.

Assume now that $X$ is reflexive and (iii) holds. By Corollary \ref{cor2}
$\big(\tfrac{1}{2}\left\Vert \cdot\right\Vert ^{2}\big)^{\ast}=\tfrac{1}%
{2}\left\Vert \cdot\right\Vert _{\ast}^{2}$ is essentially Fr\'{e}chet
differentiable, and this amounts to the Fr\'{e}chet differentiability of
$\left\Vert \cdot\right\Vert _{\ast}^{2}$ on $X^{\ast};$ using Anderson's
Theorem mentioned above we obtain that $X$ is an E-space. \hfill$\square$

\medskip

In order to obtain new variational characterizations of the closed convex sets
in an E-space let us recall that, given a lsc function $I\in F(X)$, the
problem
\[
P(I):\quad\text{minimize }I(x)\text{ for }x\in X
\]
is said to be Tykhonov well posed (TWP) if $I$ attains a strong minimum over
$X$ (i.e.\ any minimizing sequence is norm-convergent, see e.g.\ \cite{R8}).
Several characterizations of the E-spaces in terms of TWP problems have been
established (see \cite[Th.\ II.2]{R8} or \cite[Th.\ 2 p.\ 150]{R11}): we know
that the Banach space $X$ is an E-space iff for any nonempty closed convex set
$K$ in $X$ the problem $P(\left\Vert \cdot\right\Vert +\iota_{K})$ or,
equivalently, $P(\tfrac{1}{2}\left\Vert \cdot\right\Vert ^{2}+\iota_{K})$, is
TWP. By \cite[Th.\ II.2]{R8}, the Banach space $X$ is an E-space iff for any
$x^{\ast}\in X^{\ast}\setminus\{0\}$ the problem $P(x^{\ast}+\iota_{B(X)})$,
where $B(X)$ denotes the closed unit ball of $X$, is TWP. The same theorem
says that, denoting $S(X)$ the unit sphere of $X$, $X$ is an E-space iff the
problem $P(x^{\ast}+\iota_{S(X)})$ is TWP for any $x^{\ast}\in X^{\ast
}\setminus\{0\}$. Our Proposition \ref{prop4} provides another
characterization of such spaces: the Banach space $X$ is an E-space iff for
any $x^{\ast}\in X^{\ast}$ the problem $P(\tfrac{1}{2}\left\Vert
\cdot\right\Vert ^{2}+x^{\ast})$ is TWP.

To end this paper let us go back to $f$-strongly Tchebychev sets (see
Definition \ref{def4}) in the case when $f=\tfrac{1}{2}\left\Vert
\cdot\right\Vert ^{2}$. In this situation a nonempty closed $S$ in $X$ is
$f$-strongly Tchebychev iff the problem $P(\tfrac{1}{2}\left\Vert
\cdot\right\Vert ^{2}+\iota_{S}-x^{\ast})$ is TWP for every $x^{\ast}\in
X^{\ast}$. If the underlying Banach space is a Hilbert space this amounts to
say that the problem
\[
\text{minimize }\left\Vert x^{\ast}-u\right\Vert \text{ for }u\in S
\]
is TWP for any $x^{\ast}\in X^{\ast}=X.$

\begin{proposition}
\label{prop6}Let $S$ be a nonempty closed set in a Banach space $X$. Assume
that for any $x^{\ast}\in X^{\ast}$ the problem
\[
\text{minimize }\tfrac{1}{2}\left\Vert x\right\Vert ^{2}-\left\langle
x,x^{\ast}\right\rangle \text{ for }x\in S
\]
is TWP. Then $S$ is convex. If $X$ is an E-space, the converse holds.
\end{proposition}

Proof. The first part follows from the first part Theorem \ref{t2} applied to
$\tfrac{1}{2}\left\Vert \cdot\right\Vert ^{2}$. Assume now that $X$ is an
E-space. By the first part of Proposition \ref{prop5} we know that $\tfrac
{1}{2}\left\Vert \cdot\right\Vert ^{2}$ is essentially firmly
sub\-differentiable and we conclude the proof with the second part of Theorem
\ref{t2}. \hfill$\square$

\end{document}